\input amstex
\documentstyle{amsppt}
%\magnification 1200
\NoBlackBoxes

\TagsOnRight

\def\cal{\Cal}
\def\AA{{\cal A}}

\def\HH{{\cal H}}

\def\FF{{\cal F}}

\def\Z{{\Bbb Z}}

\def\R{{\Bbb R}}

\def\Q{{\Bbb Q}}

\def\n{\noindent}
\def\part{{\partial}}

\rightheadtext{Spectral invariants} \leftheadtext{
Yong-Geun Oh }

\topmatter
\title
Normalization of the Hamiltonian and the action spectrum
\endtitle
\author
Yong-Geun Oh\footnote{Partially supported by the NSF Grant \#
DMS-9971446 and by \# DMS-9729992 in the Institute for Advanced
Study, and by the grant of the Korean Young Scientist Prize
\hskip10.5cm\hfill}
\endauthor
\address
Department of Mathematics, University of Wisconsin, Madison, WI
53706, ~USA \& Korea Institute for Advanced Study, Seoul, Korea
\endaddress

\abstract In this paper, we prove that the two well-known natural
normalizations of Hamiltonian functions on the symplectic manifold
$(M,\omega)$ canonically relates the action
spectra of different normalized Hamiltonians on {\it arbitrary}
symplectic manifolds $(M,\omega)$. The natural class of normalized
Hamiltonians consists of those whose mean value is zero for the
closed manifold, and those which are compactly supported in
$\text{Int} M$ for the open manifold. We also study the effect of
the action spectrum under the $\pi_1$ of Hamiltonian diffeomorphism
group. This forms a foundational basis for our study of spectral
invariants of the Hamiltonian diffeomorphism in [Oh4].
\endabstract

\endtopmatter

\bigskip

\centerline{April, 2002}
\head
\bf \S 1. Introduction
\endhead

The study of critical values of the classical action functional
is important to study the existence question of periodic
orbits in the critical point theory initiated by Rabinowitz's
celebrated work [Ra] and in its applications to symplectic
rigidity properties (see [EH1,2], [H]). More recently, it turns out to be
useful for the study of geometry of Lagrangian submanifolds
[V], [Oh1,2] or of Hamiltonian diffeomorphism
group of symplectic manifolds [Sc], [En], [Oh3].

The most well-known difficulty that has plagued the study of {\it
family of Hamiltonian diffeomorphisms} via the mini-max theory of
the action functional is the fact that the action functional on
non-exact symplectic manifolds is not single valued on the
(contractible) loop space $\Omega_0(M)$ and has to be considered
on a $\Gamma$-covering space: Let $(\gamma,w)$ be a pair of
$\gamma \in \Omega_0(M)$ and $w$ be a disc bounding $\gamma$. We
say that $(\gamma,w)$ is {\it $\Gamma$-equivalent} to
$(\gamma,w^\prime)$ iff
$$
\omega([w'\# \overline w]) = 0 \quad \text{and }\, c_1([w'\#
\overline w]) = 0 \tag 1.1
$$
where $\overline w$ is the map with opposite orientation on the
domain and $w'\# \overline w$ is the obvious glued sphere. Here
$\Gamma$ stands for the group
$$
\Gamma = {\pi_2(M)\over \text{ker\ } (\omega|_{\pi_2(M)}) \cap
\text{ker\ } (c_1|_{\pi_2(M)})}.
$$

We denote by $[\gamma,w]$ the
$\Gamma$-equivalence class of $(\gamma,w)$ and by $\pi: \widetilde
\Omega_0(M) \to \Omega_0(M)$ the canonical projection. We also
call $\widetilde \Omega_0(M)$ the $\Gamma$-covering space of
$\Omega_0(M)$. The action functional $\AA_0: \widetilde
\Omega_0(M) \to \R$ is defined by
$$
\AA_0([\gamma,w]) = -\int w^*\omega. \tag 1.2
$$
Two $\Gamma$-equivalent pairs $(\gamma,w)$ and $(\gamma,w^\prime)$
have the same action and so the action is well-defined on
$\widetilde\Omega_0(M)$. When a periodic Hamiltonian $H:(\R/\Z) \times
M \to \R$ is given, we consider the functional $\AA_H:
\widetilde \Omega(M) \to \R$ by
$$
\AA_H([\gamma,w])= \AA_0(\gamma,w) - \int H(t, \gamma(t))dt
\tag 1.3
$$
and denote the set of critical values of $\AA_H$ by $\text{Spec}(H)$.

The mini-max theory of this action functional on the $\Gamma$-covering
space has been implicitly used in the proof of Arnold's conjecture
in the literature starting from Floer's original paper [Fl].
Recently the present author has developed this mini-max theory via
the Floer homology further and applied it to the study of geometry
of Hamiltonian diffeomorphism groups [Oh3]. (Earlier the author
[Oh1,2] and Schwarz [Sc] studied the mini-max theory via the Floer
homology for the case where the action functional is single valued
without going to the covering space. See also [FH] for earlier use
of Floer homology theory in the study of symplectic capacities).
So it appears that the non-single valuedness of the action
functional is not much of obstruction. Further applications of
this mini-max theory will be contained in [Oh4].

The less noticed and seemingly inessential nuisance in the study
of Hamiltonian diffeomorphisms in this approach is that the action
functional depends on the choice of Hamiltonians, not just on the
final diffeomorphism. Even when we fix the path of Hamiltonian
diffeomorphisms generating the given Hamiltonian diffeomorphism,
there still remains an ambiguity or freedom of choosing a
constant. The so called {\it action spectrum}, i.e., the set of
critical values of the action functional depends only on the
homotopy class of the Hamiltonian path {\it but upto overall shift
of a constant}. The necessity of proper normalization first
arose in Viterbo's construction [V], using the generating functions,
of certain spectral invariants
of Lagrangian embeddings in the contangent bundle $T^*N$
parametrized by $H^*(N)$. Viterbo [V] proved that the
difference of two such invariants is independent of the choice
of generating functions but the invariants themselves are
defined upto overall shift of a constant.
To study change of such invariants under the
Hamiltonian isotopy, one has to deal with this normalization
problem and to elliminate this ambiguity. This is particulary so
if one would like to study the effect on the invariants under
the Hamiltonian loop.

There are two natural both well-known normalization of
Hamiltonians in symplectic geometry: the first one is by
restricting to the class $C^\infty_m(M)$ of Hamiltonians
whose mean values are zero, i.e., those $h$ such that
$$
\int_M h d\mu = 0 \tag 1.4
$$
where $d\mu = {\omega^n \over n!}$ is the Liouville measure. This
normalization is natural for the case of {\it closed} $M$ , i.e.,
compact $M$ without boundary. The second one is natural for the
case {\it open} $M$. It could be non-compact or compact with
non-empty boundary for the second case. In this case, the natural
class of Hamiltonians (or Hamiltonian diffeomorphisms) are those
which have compact support in the interior of $M$, which we denote
by $C^\infty_c(M)$.

However it is not obvious what overall effect on the
action spectrum will be by restricting to these subset. For
example, the following is the obvious question to ask. We denote
$F\sim G$ iff $ \phi^1_F = \phi^1_G$ and $\{\phi_F^t\}_{0\leq
t\leq 1}$ and $\{\phi_G^t\}_{0 \leq t\leq 1}$ are isotopic.

\medskip
\n{\bf Question 1.} Let $F,\, G: M \times [0,1] \to \R$ be two
time-dependent Hamiltonians with $F \sim G$ and satisfying one of
the above normalization conditions. Does this imply
$$
\text{Spec}(F) = \text{Spec}(G) \tag 1.5
$$
as a subset of $\R$?
\medskip

This question has been studied by Schwarz [Sc] and by Polterovich
[Po] on {\it closed} manifold with normalization (1.4) for the
case of {\it symplectically aspherical} $(M,\omega)$ i.e., with
$\omega|_{\pi_2(M)} = 0$ and $c_1|_{\pi_2(M)} = 0$. In this case,
as far as the definition of the action functional is concerned,
we do not need to go to the covering space $\widetilde \Omega_0(M)$
but just work on $\Omega_0(M)$. One very interesting result proven
by Schwarz [Sc] is that there is no monodromy effect on the
spectrum  in this case and so he answered affirmatively to
the above question on $\Omega_0(M)$ without assuming
$F \sim G$ as long as $\phi^1_F = \phi^1_G$.

The first result of this paper answers the question affirmatively
for both of the above two normalizations in general in the above
setting on the covering space $\widetilde \Omega_0(M)$, which
naturally appears in the Floer theory on non-exact symplectic
manifolds.

\proclaim{Theorem I}
Let $F,\, G$ be two normalized Hamiltonians
in $\HH_m(\phi)$ on closed
$M$ (resp. in $\HH_c(\phi)$ on open $M$) with $F \sim G$.
Then
$$
\text{\rm Spec}(F) = \text{\rm Spec}(G)
\tag 1.6
$$
as a subset $\R$.
\endproclaim

This enables us to safely analyze the behavior of  homologically
essential critical values of the action functional constructed via
the Floer theory for a family of Hamiltonians functions which
appear in the study of geometry of the group of Hamiltonian
diffeomorphisms [Sc], [Oh3]. For example, combined with the fact
that $\text{Spec}(H)$ is of measure zero subset (see [Oh3]), the
following is the prototype of results which are used frequently in
the mini-max theory of action functional.

\proclaim{Corollary} Suppose that $\{ F^s\}_{s \in [0,1]}$ is a
one parameter family of normalized Hamiltonians generating the
same Hamiltonian diffeomorphism group $\phi$. If $c:[0,1] \to \R$
is a continuous function with
$$
c(s) \in \text{\rm Spec}(F^s) \subset \R \tag 1.7
$$
for all $s \in [0,1]$, then $c$ must be a constant function.
\endproclaim

One natural construction of such continuous functions is the
assignment of the homologically essential critical value
$\rho(H;a)$ for each quantum cohomology class $a \in QH^*(M,\Q)$
which are constructed in [Oh3,4] using the Floer
theory (see [Oh1,2], [Sc] for the case without quantum effect):
e.g., the assignment
$$
s \mapsto \rho(F^s;a).
$$
Theorem I will enable us to systematically study these invariants
in [Oh4].

One notable fact in both of the above normalizations is that the
class of Hamiltonians for each case is invariant under the action
by the symplectic diffeomorphism group. In the Lie algebra level,
they form a Lie sub-ideal of $\frak g = C^\infty(M)$ under the
natural Poisson bracket. More precisely, let $\frak h \subset
\frak g$ denote any of the two subset above and $\{\cdot, \cdot\}$
be the Poisson bracket on $\frak g$ associated to the symplectic
form $\omega$. Then we have
$$
\{\frak h, \frak g\} \subset \frak h. \tag 1.8
$$
It appears that this fact is a crucial link between the
normalization of Hamiltonians and the action spectrum (see \S 3
for the relevant discussion). For
example, it is not difficult to see that (1.8) will not hold true
for the kind of normalization
$$
\min H_t = 0
$$
on closed manifolds. This fact will become obvious once we have our
normalization result.

In \S 4, we analyze the action by $\pi_1(\HH am (M,\omega),id)$ on the
action spectrum. Let $h$ be an arbitrary loop in
$\HH am(M,\omega)$ based at the identity and $[h] \in \pi_1(\HH
am(M,\omega))$ be its fundamental class. To describe the result in
\S 4, let us recall some basic facts. Let
$$
\widetilde h: \widetilde \Omega_0(M) \to \widetilde \Omega_0(M)
$$
be a lifting of the action $h: \Omega_0(M) \to \Omega_0(M)$ (see [Se]).
This lifting always exists but is not unique. According to
Seidel [Se], the set of lifts forms a group
$$
\widetilde G \subset G \times \text{Homeo}(\widetilde \Omega_0(M))
$$
which is the set of pairs $(h,\widetilde h)$ such that $\widetilde
h$ is a lift of the $h$-action on $\Omega_0(M)$. And it has
one-one correspondence with the set of normalized Hamiltonian
fiber bundle with fiber $(M,\omega)$ over $S^2$ with a
$\Gamma$-equivalence class of section [Lemma 2.13,Se]. Here $G=
\Omega(\HH am(M,\omega),id)$, the space of loops based at the
identity. There exists an exact sequence
$$
0 \to \Gamma \to \widetilde G \to G \to 0.
$$
A standard calculation (Lemma 2.3 below) shows
$$
\widetilde h^*(d\AA_{H \# F}) = d\AA_F \tag 1.9
$$
and so the difference $\widetilde h^*(\AA_{H \# F}) - \AA_F$ is
constant on $\widetilde \Omega_0(M)$, which depends on $F, \, H$
and the lifting $\widetilde h$ in general. The following theorem
describes the dependence of this constant precisely.

\proclaim{Theorem II}
Let $h \in \Omega(\HH am(M,\omega))$ be a loop and let
$(h,\widetilde h)$ be any lift. Denote $[h,\widetilde h]
\in \pi_0(\widetilde G)$.  Then we have
$$
\widetilde h^*(\AA_{H\# F})=\AA_F + I_\omega([h, \widetilde h]) \tag 1.10
$$
for any $F \mapsto \phi$ and $H$ generating the loop $h$, where
$I_\omega([h, \widetilde h])$ is a constant depending only on $[h,
\widetilde h] \in \pi_0(\widetilde G)$ but independent of $F$. The
shift is induced by a homomorphism
$$
I_\omega: \pi_0(\widetilde G) \to \R
$$
that satisfies
$$
I_\omega([id, \gamma]) \in \Gamma_\omega. \tag 1.11
$$
for any lift $\gamma$ of the identity.
\endproclaim

We introduce
the {\it universal action functional}
$$
\AA: \widetilde{\HH am}(M,\omega) \times \widetilde \Omega_0(M) \to
\R
$$
by
$$
\AA(\widetilde \phi, [z,w]) := \AA_F([z,w])
$$
for $\widetilde \phi = [\phi,F]$. Theorem I implies that this is
well-defined if $F$ is normalized. We have a natural action of
$\pi_0(\widetilde G)$ on $\widetilde{\HH am}(M,\omega) \times
\widetilde \Omega_0(M)$
$$
[h,\widetilde h]: (\widetilde \phi,[z,w]) \mapsto ([h]\cdot
\widetilde \phi, \widetilde h \cdot [z,w]). \tag 1.12
$$

Evaluating $I_\omega$ at $[id,0]\in \widetilde{\HH am}(M,\omega)$,
we obtain
$$
I_\omega([id,0]) = \AA_H(\widetilde h\cdot [p, \widehat p]) = -
\int w_p^*\omega - \int_0^1 H(z_p(t),t)\, dt \tag 1.13
$$
where $[z_p,w_p] =\widetilde h\cdot [p, \widehat p]$ for any $p
\in M$.

It follows from (1.11) that in the weakly exact case these values
depend only on the loop $h$ not on $\widetilde h$. Therefore we
have a map
$$
\overline I_\omega: \pi_1(\HH am(M,\omega)) \to \R. \tag 1.14
$$
This is precisely the monodromy map considered by Schwarz
in [Sc]. In the aspherical case, i.e, for those $(M,\omega)$ with
$\omega|_{\pi_2(M)} = c_1|_{\pi_2(M)} = 0$, he proved that (1.14)
is trivial by relating the question to the geometry of Hamiltonian
fibrations [Sc].

\medskip
\n{\bf Question 2 } Is the map (1.14) trivial for the general
weakly exact case?
\medskip

We would like to thank the Institute for Advancd Study for the
excellent research environment and hospitality during our
participation of the program ``Symplectic Geometry and Holomorphic
Curves''. We also thank M. Entov and L. Polterovich for useful
e-mail communications.

\medskip

\n{\bf Convention.} \roster
\item The Hamiltonian vector field $X_f$ associated to a function
$f$ on $(M,\omega)$ is defined by $df = \omega(X_f, \cdot)$.
\item
The Poisson bracket $\{f, g\}$ is defined by $\{f,g\}:= df(X_g)$.
\item
The addition $F\# G $ and the inverse $\overline G$ on the set of
time periodic Hamiltonians $C^\infty(M \times S^1)$ are defined by
$$
\align
F\# G(x,t) & = F(x,t) + G((\phi_G^t)^{-1}(x),t) \\
\overline G(x,t) & = - G(\phi_G^t(x),t).
\endalign
$$
\endroster

\head \bf \S 2. The loop space and the action functional
\endhead
Let $(M,\omega)$ be any symplectic manifold, compact or not.
and $\Omega_0(M)$ be the set of contractible loops and
$\widetilde\Omega_0(M)$ be its the covering space mentioned before.

When a periodic Hamiltonian $H:M \times
(\R/\Z) \to \R$ is given, we consider the action functional $\AA_H:
\widetilde \Omega(M) \to \R$ by
$$
\AA_H([\gamma,w])= \AA_0(\gamma,w) - \int H(\gamma(t),t)dt
$$
We denote by $\text{Per}(H)$ the set of periodic orbits of $X_H$.
\medskip

\definition{Definition 2.1}  We define the {\it
action spectrum} of $H$, denoted as $\hbox{\rm Spec}(H) \subset
\R$, by
$$
\hbox{\rm Spec}(H) := \{\AA_H(z,w)\in \R ~|~ [z,w] \in
\widetilde\Omega_0(M), z\in \text {Per}(H) \},
$$
i.e., the set of critical values of $\AA_H: \widetilde\Omega(M)
\to \R$. For each given $z \in \text {Per}(H)$, we denote
$$
\hbox{\rm Spec}(H;z) = \{\AA_H(z,w)\in \R ~|~ (z,w) \in
\pi^{-1}(z) \}.
$$
\enddefinition

Note that $\text {Spec}(H;z)$ is a principal homogeneous space
modeled by the period group of $(M,\omega)$
$$
\Gamma_\omega := \{ \omega(A)~|~ A \in \pi_2(M)\} = \omega(\Gamma)
$$
and
$$
\hbox{\rm Spec}(H)= \cup_{z \in \text {Per}(H)}\text {Spec}(H;z).
$$
Recall that $\Gamma_\omega$ is either a discrete or a countable
dense subset of $\R$. The following was proven in [Oh3].

\proclaim\nofrills{Lemma 2.2.}~ $\hbox{\rm Spec}(H)$ is a measure
zero subset of $\R$ for any $H$.
\endproclaim

For given $\phi \in {\cal H}am(M,\omega)$, we denote by $F
\mapsto \phi$ if $\phi^1_F = \phi$, and denote
$$
\HH(\phi) = \{ F ~|~ F \mapsto \phi \}.
$$
We say that two Hamiltonians $F$ and $G$ are equivalent and denote
$F\sim G$ if they
are connected by one parameter family of Hamiltonians
$\{F^s\}_{0\leq s\leq 1}$ such that $F^s \mapsto \phi$ for all $s
\in [0,1]$. We denote by $[F]$ the equivalence class of $F$. Then
the universal covering space $\widetilde{{\cal  H}am}(M,\omega)$
of ${\cal  H }am(M,\omega)$ is realized by the set of such
equivalence classes. Note that the group $\Omega({\cal  H
}am(M,\omega))$ of based loops (at the identity)
naturally acts on the loop space $\Omega(M)$ by
$$
(h\cdot \gamma) (t) = h(t)(\gamma(t))
$$
where $h \in \Omega({\cal H}am (M,\omega))$ and $\gamma \in
\Omega(M)$. An interesting consequence of Arnold's conjecture is
that this action maps $\Omega_0(M)$ to itself (see e.g., [Lemma
2.2, Se]). Seidel [Lemma 2.4, Se] proves that this action (by a
based loop) can be lifted to $\widetilde\Omega_0(M)$.

Next if a Hamiltonian $H$ generating the loop $h$ is
given, the assignment
$$
z \mapsto h\cdot z \tag 2.1
$$
provides a natural one-one correspondence
$$
h: \text{Per}(F) \mapsto \text{Per}(H\# F)
\tag 2.2
$$
where $H\# F = H + F\circ (\phi_H^t)^{-1}$.

When the loop $h: S^1 \to \HH am(M,\omega)$ with $h(0) = id$ is
contractible, the above mentioned lifting to $\widetilde \Omega_0(M)$
of the action of $h$ on $\Omega_0(M)$ can be described explicitly:
Let $\widetilde h: D^2 \to \HH am(M,\omega)$ be a contraction of the
loop to the identity. This contraction provides a natural lifting
of the action of the loop $h$ on $\Omega_0(M)$ to
$\widetilde\Omega_0(M)$ which we define
$$
\widetilde h \cdot [\gamma,w] = [h\cdot \gamma, \widetilde h
\cdot w] \tag 2.3
$$
where $\widetilde h\cdot w$ is the natural map defined by
$$
(\widetilde h\cdot w)(y) = \widetilde h(y)(w(y)) \tag 2.4
$$
for $y \in D^2$. We call this lifting the {\it canonical lifting}
associated to the contraction $\widetilde h$. Then the assignment
$$
[z,w] \mapsto \widetilde h \cdot[z,w]
$$
provides a one-one correspondence
$$
\widetilde h: \text{Crit}(\AA_F) \to \text{Crit}(\AA_{H\# F}).
\tag 2.5
$$
Let $F,G \mapsto \phi$ and denote
$$
f_t = \phi_F^t, \, g_t = \phi_G^t,\, \text{and }\,  h_t = f_t
\circ g_t^{-1}. \tag 2.6
$$
Let $\{F^s\}_{0\leq s \leq 1} \subset \HH(\phi)$ with $F_1 =F$ and
$F_0 = G$. We denote $f^s_t = \phi_{F^s}^t$.

$\widetilde h$ provides a natural contraction of the loop $h$ to
the identity through
$$
\widetilde h: s \mapsto f^s\circ g^{-1}; \quad f^s\circ g^{-1}(t)
: = f^s_t\circ g_t^{-1}.
$$
where $\widetilde h$ can be considered as a map from $D^2$
because $h(0,s) \equiv id$. This contraction in turn provides a
canonical lifting of the action of the loop $h$, which we again
denote by $\widetilde h$, on $\Omega_0(M)$ to
$\widetilde\Omega_0(M)$ as in (2.6). For non-contractible loops
$h: S^1 \to \HH am(M,\omega)$, the lifting is related to
the notion of normalized Hamiltonian bundles (see \S 2 [Se]).

\proclaim{Lemma 2.3 } ~  Let $F, \, G$ be any
Hamiltonian with $F, G \mapsto \phi$, and $f_t, \, g_t$ and $h_t$
be as above. Suppose that $H:M \times [0,1] \to \R$ is the
Hamiltonian generating the loop $h$. We also denote the
corresponding action of $h$ on $\Omega_0(M)$ by $h$. Let
$\widetilde h$ be any lift of $h$ to $\text{Homeo}(\widetilde \Omega_0(M))$.
Then we have
$$
\widetilde h^*(d\AA_F) = d\AA_G  \tag 2.7
$$
as a one-form on $\widetilde \Omega_0(M)$.
In particular we have
$$
\AA_F \circ \widetilde h = \AA_G + C(F,G, \widetilde h) \tag 2.8
$$
where $C(F,G, \widetilde h)$ is a constant depending only on $F, \,G$.
\endproclaim
\demo{Proof} We recall that $d\AA_0$ is a pull-back $p*\alpha$ of the closed
one form $\alpha$ on $\Omega_0(M)$ by the natural projection
$p: \widetilde \Omega_0(M) \to \Omega_0(M)$ defined by
$$
\alpha(\gamma)(\xi) = \int_0^1 \omega(\dot \gamma(t), \xi(t)) \,dt
$$
for any vector $\xi \in T_\gamma\Omega_0(M)$. Therefore we first compute
$$
\align  h^*\alpha(\gamma)(\xi) & = \alpha(h\cdot \gamma)(Th(\xi))\\
&  = \int_0^1 \omega({d \over dt}(h\cdot \gamma)(t),
T_{\gamma(t)}h_t(\xi(t)))\, dt \\
& = \int_0^1 \omega(X_H(h_t(\gamma(t)) +
T_{\gamma(t)}(\dot \gamma(t)),T_{\gamma(t)}h_t(\xi(t)))\, dt \\
& = \int_0^1 \omega(X_H(h_t(\gamma(t)), T_{\gamma(t)}h_t(\xi(t)))
\, dt  \\
& \qquad + \int_0^1 \omega(T_{\gamma(t)}h_t(\dot \gamma(t)),
T_{\gamma(t)}h_t(\xi(t)))\, dt.
\endalign
$$
For the first term,  we recall that the Hamiltonian $H$ generating
$h_t$ is given by
$$
\eqalign{ H(x,t) & = F(x,t) + \overline G(f_t^{-1}(x), t) \cr & =
F(x,t)-G(g_t\circ f_t^{-1},t) = F(x,t)-G(h_t^{-1}(x),t) \cr}
$$
and so $H(h_t(x),t) = F(h_t(x),t) -G(x,t)$. Recall that
$\overline G(x,t)=-G(\phi^t(G)(x),t)$ is the Hamiltonian
generating $\phi^{-1}$. Therefore since $h_t$ is symplectic, we
have
$$
\align \omega(X_H(h_t(\gamma(t))), & T_{\gamma(t)}h_t(\xi(t)))  =
\omega
(h_t^*(X_{H_t})(\gamma(t)),\xi(t)) \\
& = \omega(X_{H_t(h_t)}(\gamma(t)),\xi(t))
= \omega(X_{(F_t(h_t)-G_t)}(\gamma(t)), \xi(t))\\
& = (d(F_t\circ h_t) - dG_t)(\xi(t))
\endalign
$$
For the second term of (2.8), we have
$$
\int_0^1 \omega(T_{\gamma(t)}h_t(\dot \gamma(t)),
T_{\gamma(t)}h_t(\xi(t)))\, dt = \int_0^1 \omega(\dot \gamma(t),
\xi(t))\, dt.
$$
Combining all these, we derive
$$
\align
\widetilde h^*d\AA_F(\gamma)(\xi) & = \widetilde h^*(
p^*\alpha)(\xi) - \int_0^1
d(F\circ h_t)(\gamma(t))(\xi(t))\, dt \\
& = \int_0^1 (d(F_t\circ h_t) - dG_t)(\gamma(t))(\xi(t)) \, dt + \int_0^1
\omega(\dot \gamma(t), \xi(t))\, dt \\
& \qquad - \int_0^1 d(F\circ h_t)(\gamma(t))(\xi(t))\, dt \\
& = \int_0^1 \omega(\dot \gamma(t), \xi(t))\, dt -\int_0^1
dG_t(\xi(t)) \, dt  = d\AA_G(\gamma)(\xi).
\endalign
$$
which proves (2.7).
\qed\enddemo

\proclaim{Corollary 2.4} Suppose $F, \, G$ are two Hamiltonians
such that $F, \, G \mapsto \phi$. Then the action spectra of
$F$ and $G$ coincide up to an overall shift by a constant
$C=C(F,G,\widetilde h)$.
\endproclaim

\head{\bf \S 3. Normalization of Hamiltonians and the action spectrum}
\endhead

Now we go back to the one-correspondence (2.5) and study its effect
on the action spectra. Naturally one hopes that $\widetilde h$
preserves the action, i.e., satisfies
$$
\AA_{H\# F}(\widetilde h\cdot [z,w]) = \AA_F([z,w])
\tag 3.1
$$
for any $[z,w] \in \text{\rm Crit}\AA_F$
for any $F$ and  for any lifting $\widetilde h$ of $h$,
{\it provided $h$ is contractible}.
The main theorem in this section is to prove that
this is indeed the case, if we use a sub-class of Hamiltonians
which are normalized accordingly depending on whether $M$ is
closed or open. On closed $M$, we consider
$$
C^\infty_m(M) := \{ f \in C^\infty(M) ~|~  \int_M f \, d\mu = 0\}
$$
and on open $M$, we consider
$$
C^\infty_c(M):= \{f \in C^\infty(M) ~|~ f \text{ has compact
support in }\, \text{Int} M \}.
$$
For given Hamiltonian diffeomorphism $\phi \in \HH am(M,\omega)$,
we denote by $\HH_m(\phi)$ (resp. $\HH_c(\phi)$) the set of $H: M \times [0,1]
\to \R$ with $H \mapsto \phi$ such that $H_t \in C^\infty_m(M)$
(resp. $C^\infty_c(M)$).

\proclaim{Theorem 3.1} Let $F,\, G$ be two normalized Hamiltonians
in $\HH_0(\phi)$ on closed
$M$ (resp. in $\HH_c(\phi)$ on open $M$) with $F \sim G$.
Then (3.1) holds and so the correspondence (2.5)
is action preserving which in turn implies
$$
\text{\rm Spec}(G) = \text{\rm Spec}(F)
\tag 3.2
$$
as a subset $\R$.
\endproclaim

Since the proof will be almost identical in both cases, we will
focus on the closed case and consider the set $C^\infty_m(M)$. In
the course of studying this case, we will also indicate the
modification we need for the open case.

This theorem in particular proves that
the spectrum $\text{Spec}(G)$ is
indeed an invariant of $\widetilde \phi \in \widetilde{\HH am}(M,\omega)$.

\definition{Definition 3.2 [Universal Action Functional]}
The {\it universal action functional}
$$
\AA: \widetilde{\HH am}(M,\omega) \times \widetilde \Omega_0(M) \to \R
$$ is defined by
$$
\AA(\widetilde \phi,[z,w]): = \AA_F([z,w])
\tag 3.3
$$
for any normalized representative $F$ with $\widetilde
\phi=[\phi,F]$.
\enddefinition

\definition{Definition 3.3 [Action Spectrum Bundle]}
For $\widetilde \phi=[\phi,G] \in \widetilde{\HH am}(M,\omega)$, we
define the {\it action spectrum} of $\widetilde \phi$ by
$$
\text{\rm Spec}(\widetilde \phi) := \text{\rm Spec}(G)
$$
for a (and so any) $G \in \HH_0(M)$ with $\widetilde \phi = [\phi,G]$.
We define the bundle of action spectrum of $(M,\omega)$ by
$$
\frak{Spec}(M,\omega) = \{ (\widetilde \phi,\AA(\widetilde \phi,
[z,w])) \mid d\AA_{\widetilde \phi}([z,w]) = 0 \} \subset
\widetilde{\HH am}(M,\omega) \times \R
$$
and denote by $\pi: \frak{Spec}(M,\omega)
\to \widetilde{\HH am}(M,\omega)$
the natural projection.
\enddefinition

It remains to prove Theorem 3.1. We first need some preparation
following Polterovich's discussion in [Proposition 6.1.C, Po] but
clarifying the role of (1.8) in the discussion. Let
$f^s_t$ be any two parameter family of Hamiltonian
diffeomorphisms. Denote by $X_{t,s}$ and $Y_{t,s}$ be the
Hamiltonian vector fields
$$
\aligned
{\part f^s_t \over \part t} \circ (f^s_t)^{-1} &= X_{t,s} \\
{\part f^s_t \over \part s} \circ (f^s_t)^{-1} &= Y_{t,s}
\endaligned
$$
Then the following is the key formula below (see [B] for its
derivation):
$$
{\part X_{t,s} \over \part s} = {\part Y_{t,s} \over \part t} + [X_{s,t},
Y_{s,t}].
\tag 3.4
$$
Now we represent the vector fields $X_{t,s}$ and $Y_{t,s}$ by the
corresponding normalized Hamiltonians $F(\cdot, t,s)$ and
$K(\cdot, t,s)$ which are uniquely determined by $X_{t,s}$ and
$Y_{t,s}$ respectively. Since (3.4) does not constrain anything on
the direction of $(t,s)$, the Hamiltonian $F$ and $K$ will satisfy
$$
{\part F \over \part s}(x,s,t) = {\part K \over \part t}(x,s,t)
-\{F,K\} (x,s,t) + c(s,t) \tag 3.5
$$
for some function $c:[0,1]^2 \to \R$ which depends only on
$(s,t)$. Here $\{F, K\}$ is the Poisson bracket i.e. $\{F,K\} =
dF(X_K)$ which satisfies the relation
$$
[X_F, X_K] = - X_{\{F,K\}}.
$$
Recall that  {\it if $F, \, K$ are normalized (in fact if just one
of them is normalized), then the Poisson bracket $\{F, K\}$ is
automatically normalized by Liouville's theorem}. Integrating
(3.5) over $M$ proves $c \equiv 0$. Therefore $F, \, K$ must
satisfy
$$
{\part F \over \part s}= {\part K \over \part t} -\{F,K\} \tag 3.6
$$
Note that  the same argument applies to prove (3.6) for the subset
$C^\infty_c(M)$ on open $M$, this time using the fact that
$$
\text{supp}\{F,K\} \subset \text{supp}F \cap \text{supp}K
$$
and considering the value at a point $p \in M \setminus
(\text{Supp}F \cup \text{Supp}K)$ in (3.5). Now we are ready to
prove Theorem 3.1.

\demo{Proof of Theorem 3.1}
Let $\FF = \{F^s\}_{s\in
[0,1]}$ be a path in $\HH_0(\phi)$ (resp. $\HH_c(\phi)$) such that
$F^0 =G$ and $F^1 = F$. Denote $h^s_t =f^s_t \circ g_t^{-1}$ and
$\widetilde h^s\cdot [z,w] = [h^s\cdot z, \widetilde h^s \cdot w]$
for any $z \in \text{Per}(G)$ as in (2.5)-(2.6). We will prove that
the function $\chi: [0,1] \to
\R$ defined by
$$
\chi(s) = A_{F^s}(\widetilde h^s\cdot [z,w])
$$
is constant.

Let $F, \, K: M \times [0,1]^2 \to \R$
be the unique normalized Hamiltonians  associated to $\{f^s_t\}$
as in Theorem 3.1. Since $f^s_1 \equiv \phi$ and $f^s_0 \equiv id$
for all $s \in [0,1]$, $K(x,1,s)$ and $K(x,0,s)$ must be constant
for each $s$ and in turn
$$
K(\cdot ,1,s) = K(\cdot,0,s) \equiv 0. \tag 3.7
$$
by the normalization condition. We now differentiate $\chi(s)$
$$
\chi'(s) = d\AA_{F^s}(\widetilde h^s\cdot [z,w])\Big({\part \over \part s}
(\widetilde h^s\cdot [z,w])\Big) -
\int_0^1 {\part F\over \part s}((h^s\cdot z)(t), t,s)\, dt.
$$
The first term vanishes since $\widetilde h^s\cdot [z,w]$ is a
critical point of $\AA_{F^s}$. For the second term
we note $z(t) = g_t(p)$ for a fixed point $p \in M$ of $\phi$
since $z$ is a periodic solution of $\dot x = X_{G_t}(x)$. Therefore
we have
$$
(h^s\cdot z)(t) = h^s_t(z(t)) = (f^s_t\circ g_t^{-1})\circ g_t(p)
= f^s_t(p).
$$
Hence we have
$$
\int_0^1 {\part F\over \part s}((h^s\cdot z)(t), t,s)\, dt =
\int_0^1 {\part F\over \part s}(f_t^s(p), t,s)\, dt.
$$
The vanishing of this integral is precisely [Lemma 12.4, Po]. For the
completeness' sake, we give its proof here. We derive from (3.2)
$$
\align
{\part F\over \part s}(f_t^s(p), t,s) & = {\part K \over
\part t}(f_t^s(p),t,s) - \{F, K\} (f^s_t(p),t,s)\\
& = {\part K \over
\part t} (f_t^s(p),t,s)) + dK(X_F)(f^s_t(p)) \\
& = {\part \over \part t}(K(f_t^s(p),t,s)).
\endalign
$$
Integrating  the
total derivative over $[0,1]$ and using (3.7), we derive
$$
\int_0^1 {\part F\over \part s}((h^s\cdot z)(t), t,s)\, dt
 = K(f_1^s(p),1,s) - K(f_0^s(p),0,s) = 0
$$
This proves that $\chi$ must be constant. Therefore (3.1) immediately
follows from $\chi(0)=0$.
(3.2) follows from the fact that both $F$ and $H\# G$
generate $g_t$ and are normalized, which in turn implies $F=H\#G$.
\qed\enddemo

\head{\bf \S 4. Action of $\pi_0(\widetilde G)$ on the action
spectrum}
\endhead

From now on, we will always assume that our Hamiltonians are
normalized.  In this section, we study effects of the action of
{\it non-contractible} loops $h$ on the action spectrum. Consider
the action spectrum bundle
$$
\pi: \frak{Spec}(M,\omega)
\to \widetilde{\HH am}(M,\omega).
$$

Let $h$ be an arbitrary loop in $\HH am(M,\omega)$ based at the
identity and $[h] \in \pi_1(\HH am(M,\omega))$ be its fundamental
class. Let
$$
\widetilde h: \widetilde \Omega_0(M) \to \widetilde \Omega_0(M)
\tag 4.1
$$
be a lifting of the canonical action $h: \Omega_0(M) \to \Omega_0(M)$
defined by
$$
(h\cdot z)(t) = h(t)(z(t)).
$$

Let $H\mapsto h$ and $(h,\widetilde h) \in \widetilde G$ be a lift
of $h$. For any $\widetilde \phi=[\phi,F]$, we consider the action
of $(h,\widetilde h)$ on $\widetilde \Omega_0(M)$,
$$
[z,w] \mapsto \widetilde h\cdot [z,w].
$$
We now study how this acts on $\text{Spec}(\widetilde \phi)$.
Since $\widetilde h^*(d\AA_{H\# F}) = d\AA_F$ by Lemma 2.3, the difference
$$
\AA_{H\# F}(\widetilde h \cdot [z,w]) - \AA_F([z,w])
$$
is independent of $[z,w]$.  We denote the common
value by $I_\omega(h,\widetilde h;F)$.

\proclaim{Lemma 4.1} The value $I_\omega(h,\widetilde h;
F)$ is independent of $F$ but depends only on $(h,\widetilde h)$.
We define
$$
I_\omega(h,\widetilde h) := I_\omega(h,\widetilde h; F)
\tag 4.2
$$
for any $F$ (e.g., $F=0$).
Then the action (4.1) induces a canonical isometry between
$\text{\rm Spec}(\widetilde \phi)$ and $\text{\rm Spec}(\widetilde
h \cdot \widetilde \phi)$ which is just the restriction of
translation by $I_\omega(h,\widetilde h)$ on $\R$.
Furthermore, $I_\omega([id,\gamma]) \equiv 0$ for any canonical lift
$\gamma$ of $id$.
\endproclaim
\demo{Proof} A straightforward computation, using the identity
$$
H\#
F(h\cdot z,t) - F(z,t) = H(h\cdot z,t),
$$
shows
$$
\aligned
I_\omega(h,\widetilde h;F) & =
\AA_{H\# F}(\widetilde h \cdot [z,w]) - \AA_{F}([z,w]) \\
& = - \int(\widetilde h\cdot w)^*\omega + \int H( h\cdot
z) - \int w^*\omega \\
& = \AA_H(\widetilde h\cdot [z,w]) - \AA_0([z,w]) = I_\omega(h,\widetilde h)
\endaligned
$$
This proves the first statement. Then the second immediately follows
from
$$
I_\omega(h,\widetilde h;F) =
\AA_{H\# F}(\widetilde h \cdot [z,w]) - \AA_{F}([z,w]).
$$
Finally $I_\omega(id,\gamma) =0$ for any canonical lift $\gamma$ of $id$
is just a restatement of Theorem 3.1. \qed\enddemo

\proclaim{Proposition 4.2} Let $\widetilde \phi \in \widetilde{\HH
am}(M,\omega)$.  $I_\omega(h,\widetilde h; \widetilde \phi)$ depends only
on $[h, \widetilde h] \in \pi_0(\widetilde G)$, i.e.,
if there is a continuous path
$\{(h^s, \widetilde h^s)\}_{s\in [0,1]}$ in $\widetilde G$, then
$$
I_\omega(h^0,\widetilde h^0) = I_\omega(h^1, \widetilde h^1). \tag 4.3
$$
And the induced map
$$
I : \widetilde{\HH am}(M,\omega) \to \R
$$
defines a group homomorphism, and
satisfies
$$
I_\omega([h,\widetilde h]) - I_\omega([h,\widetilde h']) \in
\Gamma_\omega \tag 4.4
$$
where $(h,\widetilde h)$ and $(h,\widetilde h')$ are the lifts of
the same $h$.
\endproclaim
\demo{Proof} The same proof as that of Theorem 3.1 proves (4.3).
For the homomorphism property, it is enough consider
$I_\omega(h,\widetilde h;F)$. In this case, we have
$$
\align I_\omega(h_2\cdot h_1, \widetilde{h_2\cdot h_1}); F)
& = \AA_{H_2\#H_1\# F}(\widetilde {h_2\cdot h_1} \cdot [z,w])
- \AA_F([z,w]) \\
&=(\AA_{H_2\#H_1\# F}(\widetilde {h_2\cdot h_1} \cdot [z,w])
- \AA_{H_1\# F}(\widetilde h_1\cdot [z,w])) \\
& \quad + (\AA_{H_1\# F}(\widetilde h_1 \cdot [z,w]) -
\AA_F([z,w])) \tag 4.5
\endalign
$$
for any $(h_2, \widetilde h_2)\cdot (h_1, \widetilde h_1) =
(h_2\cdot h_1, \widetilde{h_2\cdot h_1})$. The second term of (4.5)
gives rise to $I_\omega([h_1, \widetilde{ h_1}];\widetilde \phi)$. On the
other hand, we have
$$
\align \AA_{H_2\#H_1\# F} & (\widetilde {h_2\cdot h_1} \cdot
[z,w]) - \AA_{H_1\# F}(\widetilde h_1\cdot [z,w])\\
& = \AA_{H_2\#(H_1\# F)}(\widetilde h_2 \cdot (\widetilde h_1
\cdot [z,w])) - \AA_{H_1\# F}(\widetilde h_1\cdot [z,w])) =
I_\omega(h_2, \widetilde h_2)
\endalign
$$
where the last identity comes from definition of $I_\omega(h,\widetilde h;F)$
and its independence of $F$ (Lemma 4.1).

For the proof of (4.4), we just note that
$$
[h,\widetilde h'] = [id,\gamma]\cdot[h,\widetilde h]
$$
where $\gamma = \widetilde h'\circ (\widetilde h)^{-1}$.  Then it
follows from the homomorphism property of $I$ that
$$
I_\omega([h,\widetilde h'])  - I_\omega([h,\widetilde h])= I_\omega([id,\gamma]).
$$
On the other hand, we have
$$
I_\omega([id,\gamma]) = I_\omega([id,\gamma];0) = \AA_0(\gamma \cdot
[p,\widehat p]) = -\omega ([\gamma \cdot \widehat p]) \in
\Gamma_\omega
$$
where $\gamma \cdot \widehat p$ is a disc with the constant
boundary map $\widehat p$ and so defines a map from the sphere.
This finishes the proof of (4.4). \qed\enddemo

From the construction, $I_\omega([h, \widetilde h])$
is the amount of the shift in $\R$ for the isometry
$$
\widetilde h: \text{Spec}(\widetilde \phi) \to
\text{Spec}([h] \cdot \widetilde \phi); \,
[z,w] \mapsto \widetilde h\cdot [z,w]
\tag 4.6
$$
for any $\widetilde \phi \in \widetilde{\HH am}(M,\omega)$.
What we have
established by now is the following generalization of Theorem 3.1
for non-contractile loops.

\proclaim{Theorem 4.3}
Let $h \in \Omega(\HH am(M,\omega))$ be a loop with $[h,\widetilde h]
\in \pi_0(\widetilde G)$. Then the amount of shift is
independent of the lift $\widetilde h$ but only on its homotopy class
$[h,\widetilde h] \in \pi_0(\widetilde G )$.
In other words, we have
$$
\widetilde h^*(\AA_{H\# F})=\AA_F + I_\omega([h,\widetilde h]) \tag 4.7
$$
for any $F \mapsto \phi$, $H$ generating the loop $h$ and for any
lift $(h,\widetilde h)$. \endproclaim

Now let us consider the values $I_\omega([id,\gamma])$.
In the weakly exact case where $\Gamma_\omega = 0 $, (4.4) implies
that $I$ depends only on $[h] \in \pi_1(\HH am(M,\omega))$, but
not on the lift $\widetilde h \in \text{Homeo}(\widetilde
\Omega_0(M))$. Therefore it defines a map
$$
\overline I_\omega: \pi_1(\HH am(M,\omega)) \to \R
$$
defined by
$$
\overline I_\omega(a) :=I_\omega([h,\widetilde h]) = \AA_H(\widetilde h\cdot
[p,\widehat p]) = - \int w_p^*\omega - \int_0^1 H(z_p(t),t)\, dt
$$
for $a = [h]$ and $(h,\widetilde h)$ is any lifting of $h$, and
$[z_p,w_p] = \widetilde h\cdot [p, \widehat p]$. This is precisely
the monodromy map considered by Schwarz [Sc] on the weakly exact
$(M,\omega)$. He proved that this is trivial for the
symplectically aspherical case, i.e under the additional
assumption $c_1=0$. It would be interesting to see if this holds
true for arbitrary weakly exact $(M,\omega)$.

\definition{Remark 4.4} It turns out that the homomorphism
$$
I_\omega: \pi_0(\widetilde G) \to \R
$$
has the interpretation in terms of the Hamiltonian fibration:
Each Hamiltonian loop $h$ associates a Hamiltonian fibration
$\pi: P \to S^2$ with $h$ as the gluing map. Then each lift
$\widetilde h$ corresponds to a section $s: S^2 \to P$ (see [Se]).
We consider a symplectic connection of $P$ and its associated
{\it coupling form} $\omega_P$. According to Seidel [Se],
the pair $(P,s)$ is called a normalized Hamiltonian bundle.
We normalize the area of $S^2$
to be 1. Then we have the identity
$$
I_\omega([h,\widetilde h]) = \int s^* \omega_P
$$
(see [Sc]). In terms of this picture, the homomorphism property
can be understood as the gluing formula for the integral $\int s^*
\omega_P$ of the normalized Hamiltonian bundles $(P,s)$ (see [Se],
[LMP] for the discussion of the gluing).
\enddefinition

\head {\bf References}
\endhead


\widestnumber\key{FOOO} \Refs\nofrills{}

\ref\key B \by Banyaga, A. \paper Sur la structure du groupe des
diff\'eomorphismes qui pr\'eservent une forme symplectique
\jour Comm. Math. Helv. \vol 53 \yr 1978 \pages 174-227
\endref

\ref\key EH \by Ekeland, I, Hofer, H. \paper Symplectic topology
and Hamiltonian dynamics I \& II \paperinfo Math. Z.  200, (1989),
355-378, \& Math. Z. 203, (199), 553-569
\endref

\ref\key En \by Entov, M. \paper $K$-area, Hofer metric and
geometry of conjugacy classes in Lie groups \jour Invent. Math.
\yr 2001 \vol 146 \pages 93-141
\endref

\ref\key Fl \by Floer \paper Symplectic fixed points and
holomorphic spheres \jour Commun. Math. Phys. \vol 120 \pages
575-611 \yr 1989
\endref

\ref\key FH \by Floer, A., Hofer, H. \paper Symplectic homology I
\jour Math. Z. \vol 215 \pages 37-88
\endref

\ref\key H \by Hofer, H. \paper On the topological properties of
symplectic maps \jour Proc. Royal Soc. Edinburgh \vol 115 \yr 1990
\pages 25-38
\endref

\ref\key LMP \by Lalonde, F., McDuff, D., Polterovich,L.
\paper Topological rigidity of Hamiltonian loops and
quantum homology \jour Invent. Math. \vol 135 \yr 1999
\pages 369-385
\endref

\ref\key Oh1 \by Oh, Y.-G. \paper Symplectic topology as the
geometry of action functional, I \jour Jour. Differ. Geom. \vol 46
\yr 1997 \pages 499-577
\endref

\ref\key Oh2 \by Oh, Y.-G. \paper Symplectic topology as the
geometry of action functional, II \jour Commun. Anal. Geom. \vol 7
\pages 1-55 \yr 1999
\endref

\ref\key Oh3 \by Oh, Y.-G. \paper Chain level Floer theory and
Hofer's geometry of the Hamiltonian diffeomorphism group
\paperinfo preprint, April 2001,  math.SG/0104243
\endref

\ref\key Oh4 \by Oh, Y.-G. \paper Mini-max theory,
spectral invariants and geometry of the Hamiltonian diffeomorphism
group
\paperinfo preprint, June 2002
\endref

\ref\key Po \by Polterovich, L. \inbook The Geometry of the Group
of Symplectic Diffeomorphisms \publ Birkh\"auser \yr to appear
\endref

\ref\key Ra \by Rabinowitz, P. \jour Comm. Pure Appl. Math.
\vol 31 \pages 157-184 \yr 1978 \paper Periodic solutions of Hamiltonian
systems
\endref

\ref\key Sc \by Schwarz, M. \jour Pacific J. Math. \yr 2000 \vol
193 \pages 419-461 \paper On the action spectrum for closed
symplectically aspherical manifolds
\endref

\ref\key Se \by Seidel, P. \paper $\pi_1$ of symplectic
diffeomorphism groups and invertibles in quantum homology rings
\jour GAFA \yr 1997 \pages 1046-1095
\endref

\ref\key V \by Viterbo, C. \paper Symplectic topology as the
geometry of generating functions \jour Math. Ann. \vol 292 \yr
1992 \pages 685-710
\endref
\endRefs
\enddocument